\documentclass[11pt,a4paper]{amsart}
\usepackage{amssymb}
\usepackage{graphicx}

\newcommand{\eref}[1]{(\ref{#1})}

\def\cv{\ensuremath{\text {Cor}}}
\def\ld{\ensuremath{\text{LD}}}

\theoremstyle{plain}
\newtheorem{Main}{Theorem}

\newtheorem{Thm}{Theorem}

\theoremstyle{remark}
\newtheorem{Rem}[Thm]{Remark}

\def\real{{\mathbb R}}

\def\B{\mathcal B}
\def\E{\mathcal{E}}
\def\H{\mathcal{H}}
\def\O{\mathcal{O}}

\def\R{\mathcal{R}_{\epsilon,\de}}
\def\S{\mathcal{S}}

\def\F{\tilde{F}_2}

\def\const{\operatorname{const}}

\def\dist{\operatorname{dist}}

\def\id{\operatorname{id}}

\def\min{\operatorname{min}}
\def\mod{\operatorname{mod}}

\def\al{\alpha}

\def\de{\delta}

\def\vap{\varphi}
\def\th{\theta}

\def\la{\lambda}
\def\La{\Lambda}

\begin{document}

\title[Viana maps with singularities]
{Ergodic properties of Viana-like maps with singularities in the base dynamics}
\author{Jos\'e F. Alves and Daniel Schnellmann}
\address{Jos\'e Ferreira Alves, Departamento de Matem\'atica, Faculdade de Ci\^encias da Universidade do Porto,
Rua do Campo Alegre 687, 4169-007 Porto, Portugal}
\email{jfalves@fc.up.pt}
\address{Daniel Schnellmann, Ecole Normale Sup\'erieure, D\'epartment de Math\'ematiques et Applications (DMA),
45 rue d'Ulm, 75230 Paris cedex 05, France}
\email{daniel.schnellmann@ens.fr}
\date{\today}

\thanks{JFA was partially supported by Funda\c c\~ao Calouste Gulbenkian, by the European Regional Development Fund through the programme COMPETE and by the Portuguese Government through the FCT under the projects PEst-C/MAT/UI0144/2011 and PTDC/MAT/099493/2008. 
DS is supported by the Swedish Research Council.}



\begin{abstract}
We consider two examples of Viana maps for which the base dynamics has singularities (discontinuities or critical points) and show the existence of a
unique absolutely continuous invariant probability measure and related ergodic properties such as stretched exponential decay of
correlations and stretched exponential large deviations.
\end{abstract}

\maketitle

\tableofcontents

\section{Introduction}

If $\mu$ is an invariant measure for a map $F:M\to M$, its {\em basin} is the set of all points $x\in M$ such that
$$
\frac{1}{n}\sum_{i=0}^{n-1}\delta_{F^i(x)}\overset{\text{weak*}}{\longrightarrow}\mu,\qquad\text{as}\ n\to\infty.
$$
Assuming $M$ endowed with a Riemannian structure and a volume form extended to the Borel sets in $M$ (Lebesgue measure), we say that an invariant probability measure $\mu$ is a \emph{Sinai-Ruelle-Bowen (SRB) measure} if its basin has positive Lebesgue measure. It follows from Birkhoff Ergodic Theorem that ergodic absolutely continuous probability measures (acip) are necessarily SRB measures. Here, absolute continuity is always considered with respect to Lebesgue measure. We are interested in studying some statistical features of ergodic acips for certain classes of dynamical systems.

Let $\mathcal B_{1}$ and $\mathcal B_{2}$ denote Banach spaces of real valued measurable functions defined on~\( M \).
We denote the \emph{correlation} of non-zero functions $\varphi\in \mathcal B_{1}$ and  \( \psi\in \mathcal B_{2} \) with respect to a measure $\mu$ as
\[
\cv_\mu (\varphi,\psi\circ F^n):=\frac{1}{\|\varphi\|_{\mathcal B_{1}}\|\psi\|_{\mathcal B_{2}}}\left|\int \varphi\, (\psi\circ F^n)\, d\mu-\int  \varphi\, d\mu\int
\psi\, d\mu\right|.
\]
We say that we have \emph{decay
of correlations}, with respect to the measure $\mu$, for observables in $\mathcal B_1$ \emph{against}
observables in $\mathcal B_2$ if, for every $\varphi\in\mathcal B_1$ and every
$\psi\in\mathcal B_2$ we have
 $$\cv_\mu(\varphi,\psi\circ f^n)\to 0,\quad\text{ as $n\to\infty$.}$$
  Given  $\varphi\in\mathcal B_1$ and
 \( \epsilon>0 \)
we define the \emph{large deviation} of $\varphi$ at time~$n$
~as
\[
\ld_{\mu}(\varphi,\epsilon, n):=\mu\left(\left|\frac1 n \sum_{i=0}^{n-1}
\varphi\circ f^i-\int\varphi d\mu \right|>\epsilon\right).
\]
By Birkhoff's ergodic theorem  the quantity \( \ld_{\mu}(\varphi,\epsilon,n) \to 0 \), as
\( n\to \infty \), and a relevant question also in this case  is the rate of this decay.
 In our main results we shall consider $\mathcal B_1=\H_\gamma$, the space
of \emph{H\"older continuous functions} with H\"older constant $\gamma>0$. The H\"older norm of an observable $\vap\in\H_\gamma$ is given by
$$
\|\vap\|_{\H_\gamma}=\|\vap\|_\infty+\sup_{y_1\neq y_2}\frac{|\vap(y_1)-\vap(y_2)|}{\dist(y_1,y_2)^\gamma}.
$$
In such cases we shall take $\mathcal B_2 = L^\infty(\mu)$. In the proof of Theorem~\ref{t.main1} we shall also consider $\mathcal B_1$ as the space of \emph{bounded variation functions} and $\mathcal B_2$ as $L^1(\mu)$.

The purpose of this paper is to apply the theories developed in \cite{ALP}, \cite{G}, and \cite{AFLV} to two examples
of Viana maps studied in \cite{S1} and \cite{S2}, and to deduce the existence of a unique absolutely continuous invariant
probability measure and estimates for the Decay of Correlations and Large Deviations with respect to that measure. The Central Limit Theorem, the Almost Sure Invariance Principle, the Local Limit Theorem and the Berry-Esseen Theorem will also be deduced for our systems. Due to the technical nature of some of these concepts we introduce them in an appendix.

We shall consider skew-product maps similar to Viana maps. In one case with $\beta$-transformations as the base dynamics in the circle $S^1$, and in another case with a quadratic map as the base dynamics in an interval $I$.
In the following let $Q_a(x)=a-x^2$, $x\in\real$, be a Misiurewicz-Thurston quadratic map, i.e.,  the parameter $a\in(0,2)$ is chosen
such that the critical point of $Q_a$ is pre-periodic (but not periodic). The full quadratic map $Q_2$ is excluded since we
look at perturbations of the parameter~$a$. Furthermore, we assume that $Q_a$ is non-renormalizable.

\subsection{$\beta$-transformations in the base dynamics} In \cite{S1}, the map under consideration is of the form
$F_1: S^1\times\real\to S^1\times\real$:
$$
F_1(\th,x)=(\beta\th\mod1, Q_a(x)+\al\sin(2\pi\theta)),
$$
where $\beta$ is a real number greater or equal to some lower bound $\beta_a<2$ (depending on the parameter $a$ of the quadratic map $Q_a$).
This map is similar to the maps studied in \cite{V} and \cite{BST} but allowing a discontinuity in the base dynamics. If $p<0$ denotes the negative fixed point of $Q_a$, it is easy to check that there is an open interval $J\supset [p,-p]$ such that
$F_1(S^1\times J)\subset S^1\times J$,  whenever $\al>0$ is sufficiently small.

For sufficiently small $\al>0$, it is shown in \cite{S1} that the map $F_1$ is non-uniformly expanding, and furthermore  $F_1$ admits a unique acip for \emph{almost all}
$\beta\ge \beta_a$.  In this paper we improve
this result by showing that a unique acip exists, in fact, for {\em all} $\beta\ge \beta_a$. In addition, we obtain
several statistical properties for this acip, such as stretched exponential decay of correlations and
stretched exponential large deviations.

\begin{Main}
\label{t.main1}
For all small enough $\al>0$ and all $\beta>\beta_a$,
the map $F_1:S^1\times J\to S^1\times J$ admits a unique acip $\mu$ whose basin has full Lebesgue measure in $S^1\times J$. Moreover,
\begin{enumerate}
  \item
  there exist $C,\tau>0 $ such that  
$\cv_{\mu}(\varphi, \psi\circ F_1^n) \leq C e^{-\tau n^{1/3}}$ for all $\varphi\in \mathcal H_{\gamma} $ and all $\psi\in L^\infty(\mu)$;
  \item
for all $\epsilon>0$ and all $\varphi\in \mathcal H_{\gamma} $ there exist
 \(\tau'=\tau'(\tau,\varphi, \epsilon)>0 \) and
  \( C'=C'(\varphi,\epsilon)> 0 \) such that
\( LD_{\mu}(\varphi, \epsilon, n) \leq C' e^{-\tau' n^{1/7}};
 \)
 \item the Central Limit Theorem, the vector-valued Almost Sure Invariance Principle, the Local Limit Theorem and the Berry-Esseen Theorem hold for certain H\"older observables.
\end{enumerate}
\end{Main}

\subsection{Quadratic maps in the base dynamics}
Let $Q_b(\th)=b-\th^2$, $\th\in\real$ and $b\in(0,2]$, be another Misiurewicz-Thurston map and set $I=[Q_b^2(0),Q_b(0)]$.
The map studied in \cite{S2} is of the form $F_2:I\times\real\to I\times\real$:
$$
F_2(\th,x)=(Q_b^k(\th),Q_a(x)+\al s(\th)),
$$
where $k\ge1$ is an integer and $s:I\to[-1,1]$ is a coupling function which is a priori not fixed. 
Again, if $p<0$ denotes the negative fixed point of $Q_a$, it is easy to check that there is an open interval $J\supset [p,-p]$ such that
$F_2(I\times J)\subset I\times J$  whenever $\al>0$ is sufficiently small.

In \cite{S2} it is shown that there is an integer $k_0\ge1$ and a family of (non-constant)
coupling functions $s$ which are $C^2$ outside a finite number of singularities such that for each such coupling function $s$,
all $k\ge k_0$, and all sufficiently small $\al$ the map $F_2$ is non-uniformly expanding. 
In fact, the only singularities for $s$ are square root singularities.
Without loss of generality we assume that the map $Q_b$ is non-renormalizable from which follows that $Q_b^k$
has a unique acip for all $k\ge1$. (Otherwise we can restrict the map $F_2$ to a smaller region $\tilde{I}\times\real$ such that
$Q_b^k:\tilde{I}\to\tilde{I}$ and $Q_b^k$ admits a unique acip)
In this paper we will show furthermore that $F_2$ admits a unique acip with
the same statistical properties as for the map $F_1$.

\begin{Main}
\label{t.main2}
For small enough $\al>0$,
the map $F_2:I\times J\to I\times J$ admits a unique acip $\mu$ whose basin has full Lebesgue measure in $I\times J$. 
Moreover, there exists $\tau>0$ such that for all $0<\zeta<1/9$
\begin{enumerate}
\item
   there exists $C>0 $ such that  
$\cv_{\mu}(\varphi, \psi\circ F_2^n) \leq C e^{-\tau n^{\zeta}}$ for all $\varphi\in \mathcal H_{\gamma} $ and all $\psi\in L^\infty(\mu)$;
  \item
for all $\epsilon>0$ and all $\varphi\in \mathcal H_{\gamma} $ there are
 \(\tau'=\tau'(\tau,\varphi, \epsilon)>0 \) and
  \( C'=C'(\varphi,\epsilon)> 0 \) such that
\( LD_{\mu}(\varphi, \epsilon, n) \leq C' e^{-\tau' n^{\zeta'}}
 \), where $\zeta'=\zeta/(\zeta+2)$;
\item the Central Limit Theorem, the vector-valued Almost Sure Invariance Principle, the Local Limit Theorem and the Berry-Esseen Theorem hold for certain H\"older observables.
\end{enumerate}

\end{Main}

\subsection{Strategy}To prove the first two items of Theorem~\ref{t.main1} and Theorem~\ref{t.main2} we will apply the result in \cite{G} which shows the existence of a \emph{Young tower} or \emph{Gibbs-Markov structure} for
the maps $F_1$ and $F_2$, with stretched exponential tail estimates for the expansion and slow recurrence tails. These objects will be defined precisely in Section~\ref{se.NUE}, and in Section~\ref{Bounds} we obtain stretched exponential bounds on these tails. From the existence of such a tower it then follows the
decay of correlations conclusions as stated in the two theorems. The conclusions on the large deviations are an immediate consequence of \cite[Theorem~D(2)]{AFLV}. Finally, in Section~\ref{se.transitivity} we obtain the topological transitivity of the maps, which assures the uniqueness of the acip in both cases.

The third items of Theorem~\ref{t.main1} and Theorem~\ref{t.main2} follow as a direct application of Corollaries~B1-B4 in \cite{AFLV}.

\section{Non-uniform expansion and slow recurrence}\label{se.NUE}
\label{s.preliminaries}
Let $M$ be equal to $M_1=S^1\times J$ or $M_2=I\times J$ and  $F:M\to M$ be equal to $F_1$ or $F_2$, respectively. Let $\S$ be some closed set of zero Lebesgue measure of singularities/criticalities
such that $F:M\setminus\S\to M$ is a $C^2$ local diffeomorphism.
We say that $F$ is {\em non-degenerate} close to $\S$ if there are constants $B>1$ and $\xi>0$ such that the following three conditions hold.
For all $y\in M\setminus\S$ and $v\in T_yM\setminus\{0\}$, we have
\begin{itemize}
 \item[(S1)]
\hspace{.1cm}$\displaystyle{\frac{1}{B}\dist(y,\S)^{\xi}\leq
\frac{\|DF(y)v\|}{\|v\|}\leq B\dist(y,\S)^{-\xi}}$;
\end{itemize}
and for every $y_1,y_2\in M\setminus\S$ with $\dist(y_1,y_2)<\dist(y_1,\S)/2$, we have
\begin{itemize}
 \item[(S2)]
\hspace{.1cm}$\displaystyle{\left|\log\|DF(y_1)^{-1}\|-
\log\|DF(y_2)^{-1}\|\:\right|\leq
B\frac{\dist(y_1,y_2)}{\dist(y_1,\S)^{\xi}}}$;
 \item[(S3)]
\hspace{.1cm}$\displaystyle{\left|\log|\det DF(y_1)^{-1}|-
\log|\det DF(y_2)^{-1}|\:\right|\leq
B\frac{\dist(y_1,y_2)}{\dist(y_1,\S)^{\xi}}}$.
 \end{itemize}
The critical or singular set $\S$ for the map $F_1$ is the set 
$$
\big\{\{0\}\times J\big\}\cup\big\{S^1\times\{0\}\big\},
$$
and the singular set $\S$ for the map 
$F_2$ is the set 
$$
\Big\{\bigcup_{1\le i\le m} \{b_i\}\times J\Big\}\cup\big\{I\times\{0\}\big\}, 
$$
where the points $\{b_1,...,b_m\}$ consist of
the critical points of $Q_b^k$ and the points where the coupling function $s$ is not $C^2$.
It is straightforward to check that the map $F_1:M_1\to M_1$ 
satisfies the non-degeneracy conditions (S1)--(S3), where the constant $\xi$ can
be chosen equal to $1$. Regarding the map $F_2:M_2\to M_2$, since the coupling function $s$ has only square root singularities, 
one easily checks that the non-degeneracy conditions (S1)--(S3) hold for $F_2$ with $\xi=2$.

The main result in \cite{S1} and \cite{S2}, respectively, shows that the map $F$ is {\em non-uniformly expanding}, 
i.e., there is some constant $c>0$ such that for Lebesgue almost every $y\in M$
 \begin{equation} \label{eq.nonuniform}
\liminf_{n\to +\infty}\frac{1}{n}
\sum_{j=0}^{n-1}\log\|DF(F^j(y))^{-1}\|^{-1}\geq c>0.
 \end{equation}
This implies that the {\em expansion time} function
$$
\E(y)=\min\left\{N\ge1\ :\ \frac{1}{n}\sum_{j=0}^{n-1}\log\|DF(F^j(y))^{-1}\|^{-1}\geq c/2,\ \text{for all }n\ge N\right\},
$$
is defined and finite for Lebesgue almost every $y\in M$. Given $\de>0$ we define the $\de$-{\em truncated distance}
from $y\in M$ to $\S$ as $\dist_\de(y,\S)=\dist(y,\S)$ if $\dist(y,\S)\le\de$ and $\dist_\de(y,\S)=1$ otherwise.
In the next section we will see that $F$ has {\em slow recurrence to the critical set} $\S$, i.e., given
any $\epsilon>0$ there is $\de>0$ such that
\begin{equation}
\label{eq.limsup1}
\limsup_{n\to +\infty}\frac{1}{n}
\sum_{j=0}^{n-1}-\log\dist_\de(F^j(y),\S)\leq \epsilon
\end{equation}
for Lebesgue almost every $y\in M$. It follows that the
{\em recurrence time} function
\begin{equation}
\label{eq.r}
\R(y)=\min\left\{N\ge1\ :\ \frac{1}{n}\sum_{j=0}^{n-1}-\log\dist_\de(F^j(y),\S)\leq2\epsilon,\ \text{for all }n\ge N\right\},
\end{equation}
is defined and finite for a.e. $y\in M$. 

According to the results in \cite{G}, in order to prove Theorem~\ref{t.main1} and Theorem~\ref{t.main2}, 
it is left to show that all the iterates of the map $F$ are
topologically transitive on the attractor $\Lambda=\bigcap_{n\ge0}F^n(M)$
and that there exist constants $\tau,\zeta>0$ such that for any $\epsilon>0$ 
there is $\delta>0$ such that
\begin{equation}
\label{eq.tail}
|\{y\in M\ :\ \E(y)>n\ \text{or}\ \R(y)>n\}|\le\O(e^{-\tau n^\zeta}),
\end{equation}
where $|\,.\,|$ stands for Lebesgue measure. By the technique in \cite{G}, if \eref{eq.tail} is satisfied for some constants $\tau,\zeta>0$ then the same 
constants appear also in the decay of correlations.
%

\section{Stretched exponential bounds}\label{Bounds}
The main part in the proof of the theorems is to establish the stretched exponential bound in \eref{eq.tail}.
We divide the singular sets of $F_1$ and $F_2$ into two parts. One part will contain the singularities for which it is enough to
study the base dynamics, and the other part contains the critical points due to the quadratic map $Q_a$. More precisely, when considering $F_1$
let $\S_h=\{0\}\times J$ and $\S_v=S^1\times\{0\}$ (the indices $h$ and $v$ stand for {\em horizontal} and {\em vertical}, respectively). When
considering $F_2$ let $\S_h=\bigcup_{1\le i\le m} \{b_i\}\times J$ and $\S_v=I\times\{0\}$.

Let $\mathcal{R}_{\epsilon,\de,h}$ and $\mathcal{R}_{\epsilon,\de,v}$ be defined in the same way as the set $\R$ (see \eref{eq.r})
but with $\S$ in its definition replaced by $\S_h$ and $\S_v$, respectively.
Obviously, we have
\begin{multline*}
\{y\in M\ :\ \E(y)>n\ \text{or}\ \R(y)>n\}\\
\subset\{y\in M\ :\ \mathcal{R}_{\epsilon,\de,h}(y)>n\}\cup\{y\in M\ :\ \E(y)>n\ \text{or}\ \mathcal{R}_{\epsilon,\de,v}(y)>n\}.
\end{multline*}
Hence,
in order to show \eref{eq.tail} it is sufficient to show that there exist constants  $\tau,\zeta>0$ such that for any $\epsilon>0$ there is $\delta>0$
such that
\begin{equation}
\label{eq.horizontal}
|\{y\in M\ :\ \mathcal{R}_{\epsilon,\de,h}(y)>n\}|\le\O(e^{-\tau n^\zeta}),
\end{equation}
and
\begin{equation}
\label{eq.vertical}
|\{y\in M\ :\ \E(y)>n\ \text{or}\ \mathcal{R}_{\epsilon,\de,v}(y)>n\}|\le\O(e^{-\tau n^\zeta}).
\end{equation}

\subsection{Bounds for the fiber maps}
\label{ss.vertical}
The main calculations here are done in \cite{S1} and \cite{S2} where the positivity of the Lyapunov exponents is shown.
We can essentially follow Section~6.2.1 in \cite{AA} which establishes tail estimates of the expansion
and recurrence time function for the maps studied in \cite{V} (also in their case, the essential part of the argument is done in
the proof of positive Lyapunov exponents, see \cite{V}).

We will treat the maps $F_1$ and $F_2$ simultaneously. In order to apply the results in \cite{S2}, we first have to
conjugate the function $F_2$ to a function denoted by $\F$. The conjugation function $\Phi:I\times J\to[-1,1]\times J$ is of the form
$\Phi(\th,x)=(\varphi(\th),x)$, $(\th,x)\in I\times J$, where $\varphi:I\to[-1,1]$ is analytic outside a finite number of singularities.
$\varphi$ is obtained by integrating the density of the acip for $Q_b$ from which follows that the singularities of $\varphi$ are of square root type.
The conjugation function $\Phi$ is explained in detail in \cite[p. 2684]{S2}.
The conjugated map $\F=\Phi\circ F_2\circ\Phi^{-1}$ has the form $\F:[-1,1]\times J\to[-1,1]\times J$:
$$
\F(\th,x)=(g(\th),Q_a(x)+\al h(\th)),
$$
where $g=\vap\circ Q_b^k\circ\vap^{-1}:[-1,1]\to[-1,1]$ is analytic and uniformly expanding outside a finite set of singularities and $h:[-1,1]\to[-1,1]$ is $C^2$
(extendable to a neighborhood of $[-1,1]$) with first derivative bounded away from $0$. In the setting of the map $F_1$ let
the base dynamics $\th\mapsto d\th\mod1$ also be denoted by $g$. Depending on the context, 
in the following let $M$ denote either $M_1$ or $[-1,1]\times J$, and the map $F$ stands either for $F_1$ or $\F$.
We define inductively $f_n(\th,x)$, $(\th,x)\in M$. 
$f_1(\th,x)$ is equal to $Q_a(x)+\al\sin(2\pi\th)$ for $F_1$ and equal to $Q_a(x)+\al h(\th)$ for $\F$.
For $n\ge 2$, $f_n$ is defined by the equation $F^n(\th,x)=(g^n(\th),f_n(\th,x))$.
In order to get the bound \eref{eq.vertical} of the tail of the expansion and recurrence time function, 
we have to study the returns of $f_n(\th,x)$ to $0$. 

Henceforth, we consider only points $(\th,x)\in M$ whose orbit does not hit the critical set $\S_v$. This is
no restriction since the  set of those points has full Lebesgue measure.
For $r\ge0$, set
$$
 J(r)=\{x\in I\ :\ |x|\le\sqrt\al e^{-r}\},
$$
and for each integer $j\geq 0$, we define 
$$
 r_j(\th,x)=\min\left\{r\ge0\ :\ f_j(\th,x)\in J(r)\right\}.
$$
In \cite{S1} and \cite{S2}, for some given constant $0<\kappa<1/4$, one considers
$$
 G=\left\{0\le j<n\ :\ r_j(\th,x)\ge\bigg(\frac12-2\kappa\bigg)\log\frac1\al\right\}.
$$
Fix some integer $n\ge1$ sufficiently large (only depending on $\al>0$). From the estimates in \cite[equation (14)]{S1}
and \cite[equation (23)]{S2}, we deduce that if we take
$$
 B_2(n)=\left\{(\th,x)\in M\ :\ \:\mbox{there is $1\le j<n$ with  $f_j(\th,x)\in J(\sqrt n)$\:}\right\},
$$
then there is a constant $\tau_2>0$ such that
\begin{equation}
 \nonumber
 |B_2(n)|\le\mbox{const}\,e^{-\tau_2\sqrt n}.
\end{equation}
Furthermore, there exists a constant $c>0$ (only depending on the quadratic map $Q_a$, and not on $\al$) such that
\begin{equation}
 \label{eq.principal}
 \log|\partial_xf_n(\th,x)|\ge
 cn-\sum_{j\in G}r_j(\th,x),\quad\mbox{for}\ (\th,x)\notin B_2(n),
\end{equation}
see \cite[equation (15)]{S1}, \cite[equation (24)]{S2}, and \cite[p 75 \& 76]{V}. Let
$$
 B_1(n)=\left\{(\th,x)\in M\ :\ \sum_{j\in G}r_j(\th,x)\ge \frac c2n\right\}.
$$
It is shown in \cite[equation (17)]{S1} and \cite[equation (25)]{S2} that there is $\tau_1>0$ such that
\begin{equation}
 \nonumber
 |B_1(n)|\le\mbox{const}\,e^{-\tau_1\sqrt n}.
\end{equation}
Since the base dynamics of $F$ is uniformly expanding,
we obtain immediately that
\begin{equation}
\label{eq.eee}
|\{y\in M\ :\ \E(y)>n\}|\le|B_1(n)\cup B_2(n)|\le\O(e^{-\tau\sqrt{n}}),
\end{equation}
where $\tau=\min\{\tau_1,\tau_2\}$.  Note that while the base dynamics $g$ of $\F$ is uniformly expanding the base
dynamics $Q_b^k$ of $F_2$ is not. However $Q_b^{nk}=\vap^{-1}\circ g^n\circ\vap$
and, by the properties of the density of the acip for $Q_b$ (see, e.g., \cite{S2}),
it follows that the derivative of $\vap$ is bounded away from zero (on the support of the acip) and the derivative of
$\vap^{-1}$ is strictly positive outside a finite number of critical points of order $2$.
Hence, there exists $\la>1$ such that $|D_\th Q_b^{nk}(\th)|\ge\la^n$ for all $\th$ outside an exceptional set whose
size is decreasing exponentially in $n$. It follows that the tail estimate \eref{eq.eee} of the expansion time function 
does not only hold for the maps $F_1$ and $\F$ but also for the map $F_2$.

From the arguments in \cite{S1}, \cite{S2}, and \cite{V} it is obvious that the constant $c$ in the definition of $B_1(n)$ can be
chosen arbitrarily small. Observe that in the set $B_1(n)$ we are only concerned about the returns of $f_n(\th,x)$ to the
set $J((1/2-2\kappa)\log(1/\al))$. Hence, setting
$$
\de=\frac{|J((1/2-2\kappa)\log(1/\al))|}{2}=\al^{1-2\kappa},
$$
and writing $\epsilon$ instead of $c/2$, we obtain
\begin{equation}
 \nonumber
 \sum_{j=0}^{n-1}-\log\dist_\delta(F^j(\th,x),\S_v)=\sum_{j\in G}r_j(\th,x)\le\epsilon n,
\end{equation}
for all $(\th,x)\notin B_1(n)\cup B_2(n)$.
Considering the map $F_2$ this implies that
$$
\sum_{j=0}^{n-1}-\log\dist_\delta(F_2^j(\th,x),\S_v)\le\epsilon n,
$$
for all $(\th,x)\notin\Phi^{-1}(B_1(n)\cup B_2(n))$, where $\Phi=(\varphi,\id)$ is the conjugating function described above.
Since the derivative of $\varphi^{-1}$ is bounded from above (see, e.g., \cite{S2}), we obtain
$$
|\Phi^{-1}(B_1(n)\cup B_2(n))|\le\|D\vap^{-1}\|_\infty|B_1(n)\cup B_2(n)|\le\mbox{const}\,e^{-\tau\sqrt{n}}.
$$
We conclude that for the maps $F_1$ and $F_2$ we have
$$
|\{y\in M\ :\ \mathcal{R}_{\epsilon,\de,v}(y)>n\}|\le\O(e^{-\tau\sqrt{n}}).
$$
Altogether we proved for $F_1$ and $F_2$ the stretched exponential bounds required in \eref{eq.vertical} where the constant
$\zeta$ can be taken equally to $1/2$. 

\subsection{Bounds for the base dynamics}
\label{ss.horizontal}
Note that to prove the decay on $|\{y\in M\ :\ \mathcal{R}_{\epsilon,\de,h}(y)>n\}|$ we only have to consider the base dynamics.
For the sake of notation we make this more precise.
Consider the projection of $\S_h$ to the first coordinate. We denote this projection again by $\S_h$, i.e., for the base
dynamics $g_1$ of the map $F_1$ the critical set $\S_h$ is equal to $0\in S^1$ and for the base dynamics $g_2$ of the
map $F_2$ the critical set $\S_h$ is equal to $\{b_1,...,b_m\}\subset I$.
For $y=(\th,x)\in M_i$, $i=1,2$, we have
\begin{multline*}
\mathcal{R}_{\epsilon,\de,h}(y)=\mathcal{R}_{\epsilon,\de,h}(\th)\\
:=\min\left\{N\ge1\ :\ \frac{1}{n}\sum_{j=0}^{n-1}-\log\dist_\de(g_i^j(\th),\S_h)\leq2\epsilon,\ \text{for all }n\ge N\right\},
\end{multline*}
where $\dist_\de$ is defined as above but restricted to $S^1$ or $I$, respectively.
It follows that $|\{y\in M\ :\ \mathcal{R}_{\epsilon,\de,h}(y)>n\}|$ is equal to
$|\{\th\in S^1\ :\ \mathcal{R}_{\epsilon,\de,h}(\th)>n\}||J|$ or $|\{\th\in I\ :\ \mathcal{R}_{\epsilon,\de,h}(\th)>n\}||J|$, respectively.

To establish the desired tail estimates of the recurrence time function for the base dynamics we follow the strategy of \cite[Theorem~4.2]{AFLV}.
We begin by introducing some auxiliary functions. For $\de>0$, let
\begin{equation*}
\phi(\th)=\begin{cases}
-\log\dist(\th,\S_h)&\text{ if $\dist(\th,\S_h)<\delta$ },\\
\frac{\log\delta}{\delta}(\dist(\th,\S_h)-2\delta)&\text{ if $\delta\leq \dist(\th,\S_h)<2\delta$ },\\
0&\text{ if $\dist(\th,\S_h)\geq2\delta$},
\end{cases}
\end{equation*}
where $\th$ is in $S^1$ or $I$, respectively.
Observe that $\phi$ has discontinuities at the singular set $\S_h$.
Let $\nu$ denote the unique acip for $g_1$ or $g_2$, respectively.
We can choose \( \delta>0 \) sufficiently small so that
\begin{equation*} 
    \limsup_{n\to+\infty}
\frac{1}{n} \sum_{j=0}^{n-1}- \log \dist_{\delta}(g_i^{j}(\th), \S_h)
\leq     \lim_{n\to+\infty}
\frac{1}{n} \sum_{j=0}^{n-1} \phi(g_i^{j}(\th))
=
\int \phi  d\nu \leq \epsilon,
\end{equation*}
$i=1,2$.
For all \( k >0  \)  we let
\[ A_{k} := \{\th: \phi(\th) \geq k\}
\]
and  define
\[
\phi_{k}(\th):= \begin{cases}
k, &\text{ if } \th\in A_{k};\\
\phi(\th), &\text{ otherwise}.
\end{cases}
\]
The functions $\phi_k$ and the sets $A_k$ correspond to the functions $\phi_{2,k}$ and $A_{2,k}$ in \cite[Section 5]{AFLV}, respectively.

\subsubsection{$\beta$-transformations}
We consider first the setting in the case
of the map $F_1$.
Let $\B$ denote the space of functions $\vap$ on $S^1$ with bounded variation and set
$$\|\vap\|_\B:=V_{S^1}\vap+\|\vap\|_{L^1(m)},$$ where $m$ denotes the Lebesgue measure on $S^1$.
By \cite[Appendix~C.1 and Corollary~H~(2)]{AFLV} we get that for all $\vap\in\B$ and for every $\epsilon>0$
there exists $\tau(\vap,\epsilon)>0$ and $C(\vap,\epsilon)$ such that
\begin{equation}
\label{eq.ld}
LD_\nu(\vap,\epsilon,n)=\nu\left(\left|\frac{1}{n}\sum_{i=0}^{n-1}\vap\circ g_1^i-\int\vap d\nu\right|>\epsilon\right)\le C(\vap,\epsilon)e^{-\tau(\vap,\epsilon)n}.
\end{equation}
Furthermore, by \cite[Proposition~2.5 and Lemma~2.6]{AFLV}, we get
$$
\tau(\vap,\epsilon)\ge\epsilon^2(8(\|\vap\|_\infty+C'\|\vap\|_\B)^2)^{-1}.
$$
The constant $C'$ is equal to $2\sum_{i\ge0}\xi(i)$, where $\xi(i)$ is an upper bound for the decay of
correlation for observables in $\B$ against $L^1(\nu)$. By \cite[Appendix~C.1 and Corollary~H~(1)]{AFLV}, this decay is exponential
and, hence, $C'$ is finite.
Regarding the constant $C(\vap,\epsilon)$, we derive from the proof of \cite[Proposition~2.5]{AFLV} that
$C(\vap,\epsilon)\le2e^{\epsilon(4\|\vap\|_\infty)^{-1}}$.
Since the function $\phi$ is not of bounded variation we cannot apply \eref{eq.ld} directly to $\phi$.
However, the functions $\phi_k$ are of bounded variation and,
according to \cite[equation~5.1]{AFLV}, we have
$$
LD_\nu(\phi,2\epsilon,n)\le LD_\nu(\phi_k,\epsilon,n)+n\nu(A_k).
$$
Since the density of $\nu$ is bounded from above, this immediately implies $\nu(A_k)\le\const|A_k|\le\const e^{-k}$.
Altogether we obtain,
$$
LD_\nu(\phi,2\epsilon,n)\le2e^{\epsilon(4\|\phi_k\|_\infty)^{-1}}e^{-\epsilon^2(8(\|\phi_k\|_\infty+C'\|\phi_k\|_\B)^2)^{-1}n}+\const ne^{-k}.
$$
Observe that $\|\phi_k\|_\infty=k$, $V_{S^1}\phi_k=2k$ and $\|\phi_k\|_{L^1(m)}$ is bounded from above by a constant independent on $k$.
We derive that there is a constant $C$ independent on $k$ and $\epsilon$ such that
$$
LD_\nu(\phi,2\epsilon,n)\le C(e^{-\epsilon^2C^{-1}k^{-2}n}+ne^{-k}).
$$
Choosing $k=n^{1/3}$, we get $LD_\nu(\phi,2\epsilon,n)\le\O(e^{-\tau n^{1/3}})$ for some constant $\tau=\tau(\epsilon)>0$.
Since $-\log\dist_\de\le\phi$, the density of $\nu$ is bounded away from zero (see, e.g. \cite[Appendix~C.1]{AFLV}), and
$\int\phi d\nu\le\epsilon$, we finally obtain
$$
|\{\th\in S^1\ :\ \mathcal{R}_{\epsilon,\de,h}(\th)>n\}|\le\const LD_\nu(\phi,2\epsilon,n)\le\O(e^{-\tau n^{1/3}}),
$$
which implies \eref{eq.horizontal} where $\zeta=1/3$.

\subsubsection{Quadratic maps} It is only left to consider the case of the map $F_2$. By, e.g., \cite{KN}, $(g_2,\nu)$ has exponential decay of correlations
for functions of bounded variation only against $L^p(\nu)$, $p>2$. Thus, we cannot apply the argument above for the map $g_1$ which
gives a sharper result. However, by \cite{G}, it follows that $(g_2,\nu)$ has exponential decay of correlations for H\"older
observables against $L^{\infty}(\nu)$. In order to apply \cite{G}, we need the existence of a tower with exponentially small tales. But
this follows from \cite{Y}. Moreover, by \cite{S2}, it follows that the density of $\nu$ is uniformly bounded away from $0$ on
its support. Thus, \cite[Proposition~4.1]{AFLV} implies that there exists $\tau>0$ such that for any
$0<\zeta<1/9$ and any $\epsilon>0$ sufficiently small one has
$$LD_\nu(\phi,2\epsilon,n)\le\O(e^{-\tau n^{\zeta}}).$$
This concludes the proof of \eref{eq.horizontal}
in the case of the map $F_2$.

\begin{Rem}
For $(g_2,\nu)$, \cite{KN} and \cite{MN} show exponential large deviation estimates for observables of bounded variation and for H\"older observables,
respectively. Thus, regarding the argument above for the map $g_1$ one might expect to get in \eref{eq.horizontal} a constant $\zeta$ close to 1/3.
However, the constants in \cite{KN} and \cite{MN} for the exponential large deviation are not as explicit as in \eref{eq.ld},
which makes it difficult to apply their results to our setting.
\end{Rem}

\section{Topological transitivity}\label{se.transitivity}
Denote by $\La$ the attractor $\bigcap_{n\ge0}F^n(M)$. We say that $F$ is {\em topologically transitive} on
the attractor $\Lambda$ if, for every non-empty open subsets $U$ and $V$ of $\La$, there exists $n$ such that $F^{-n}(U)\cap V$
contains a non-empty open set.
For the maps $F_i$, $i=1,2$, we have, by the same argument as in \cite[Lemma~6.1]{AV}, that its
attractor $\La_i$ coincides with $F_i^2(M_i)$. In fact, for the later use we note that the argument in \cite{AV} shows even that if
$D$ is an interval with its boundary
points sufficiently close to $Q_a^2(0)$ and $Q_a(0)$, respectively, then $F_1^2(S^1\times D)=\La_1$ and $F_2^2(I\times D)=\La_2$.

The essential part here is done in \cite{Alv}.
By Sections~\ref{ss.vertical} and \ref{ss.horizontal} we know that $F_i:M_i\to M_i$, $i=1,2$, is non-uniformly expanding and
slowly recurrent to the critical set.
Hence, we can apply Lemma~4.3 in \cite{Alv} and we get that there is a constant $\delta>0$ only dependent on the constant $c$ from the
non-uniform expansion (see equation~\eref{eq.principal}) and on the constant $\beta$ from the non-degeneracy condition such that the
following holds. For every $\epsilon>0$ there exists $n_1=n_1(\epsilon)>0$ such that for any ball $B\subset M_i$ of radius $\epsilon$ there is
an integer $n\le n_1$ such that $F_i^n(B)$ contains a ball of radius $\delta$. Recall that the constant $c$ and, thus, also the constant $\de$
do not depend on $\al$.
The following argument is similar to that in \cite[p. 29]{AV}. Recall that we defined $I=[Q_b^2(0),Q_b(0)]$.
Since $Q_a$ and $Q_b$ are non-renormalizable, it follows that the supports of the acip's
for $Q_a$ and $Q_b$ are equal to $[Q_a^2(0),Q_a(0)]$ and $I$, respectively. Since the critical points
of $Q_a$ and $Q_b$ are eventually mapped into repelling periodic points and since $Q_a$ and $Q_b$ are conjugated to uniformly
expanding maps (see, e.g., \cite[Proposition 2.2]{S2}), it follows that there is an integer $n_2=n_2(\de)>0$ such that if
$V\subset[Q_a^2(0),Q_a(0)]$ and $V'\subset I$ are intervals of length $\de$ then $Q_a^{n_2}(V)=[Q_a^2(0),Q_a(0)]$ and
$Q_b^{n_2}(V')=I$. Recall that if $D$ is an interval with its boundary
points sufficiently close to $Q_a^2(0)$ and $Q_a(0)$, respectively, then $F_1^2(S^1\times D)=\La_1$ and $F_2^2(I\times D)=\La_2$.
Since $F_1$ and $F_2$ depend continuously on $\al$, it follows that if $\th\in S^1$, $\th'\in I$, and  $V,V'$ are
intervals of length $\de$ satisfying $\th\times V\subset\La_1$ and $\th'\times V'\subset\La_2$, then for $\al$ sufficiently small we have
$$
F_1^{n_1+2}(\th\times V)=\{g^{n_1+2}(\th)\}\times\real\cap\La_1,
$$
and
$$
F_2^{n_1+2}(\th'\times V')=\{Q_b^{(n_1+2)k}(\th')\}\times\real\cap\La_2.
$$
Altogether, we derive that for each $\epsilon>0$ there is an integer $n_0=n_0(\epsilon)$, such that if $B\subset\La_i$, $i=1,2$, is a ball of radius $\epsilon$
then $F_i^{n_0}(B)=\La_i$. Thus, we conclude that $F_i$ is topologically transitive on $\La_i$.
Obviously, this argument works also for arbitrary iterates of $F_i$.

\appendix
\section{Limit Theorems}\label{appendix.further}
Here we define the statistical properties of dynamical systems which are mentioned in the last items of Theorems~\ref{t.main1} and \ref{t.main2}.

\subsection{Central Limit Theorem}%
%
%
Let $\varphi\in\mathcal H_\gamma$ be such that $\int \varphi d\mu=0$. Then
\begin{equation}\label{eq.sigma}
\sigma^2=\lim_{n\to\infty}\frac1n \int \left(\sum_{i=0}^{n-1} \varphi\circ F^i\right)^2d\mu\geq0
\end{equation}
is well defined. We say the Central Limit Theorem holds for $\varphi$ if for all $a\in\mathbb R$
\[
\mu\left(\left\{x: \frac1{\sqrt n} \sum_{i=0}^{n-1} \varphi\circ F^i(x)\leq a\right\}\right)\rightarrow \int_{-\infty}^a \frac1{\sigma\sqrt{2\pi}}\text e^{-\frac{x^2}{2\sigma^2}}dx, \text{ as $n\to\infty$},
\]
whenever $\sigma^2>0$.
Additionally, $\sigma^2=0$ if and only if $\varphi$ is a \emph{coboundary}  (\( \varphi\neq \psi\circ F - \psi \) for any \(
\psi \in L^2\).

%

\subsection{Local Limit Theorem}

\label{aperiodic} A function $\varphi:M\to\mathbb R$ is said to be \emph{periodic} if there exist $\rho\in\mathbb R$, a measuruble $\psi:M\to\mathbb R$, $\lambda>0$ and
$q:M\to\mathbb Z$, such that $$\varphi=\rho+\psi-\psi\circ F+\lambda q$$ almost everywhere. Otherwise, it is said to be \emph{aperiodic}.
%
%

Let
$\varphi \in \mathcal H_{\gamma}$ be such that $\int \varphi d\mu=0$ and $\sigma^2$ be as in \eqref{eq.sigma}. Assume that  $\varphi$ aperiodic (which implies that $\sigma^2>0$). We say that the Local Limit Theorem holds for  $\varphi$ if for any bounded interval $J\subset \mathbb R$,  for any real sequence $\{k_n\}_{n\in\mathbb N}$ with $k_n/n \to \kappa \in\mathbb R$, for any $u\in\mathcal H_{\gamma} $, for any measurable $v:M\to\R$ we have
\[
\sqrt n \mu\left(\left\{x\in M:\; \sum_{i=0}^{n-1} \varphi\circ F^i(x)\in J+k_n+u(x)+v(F^nx)\right\}\right)\to m(J) \frac{\mbox{e}^{-\frac{\kappa^2}{2\sigma^2}}}{\sigma\sqrt{2\pi}}.
\]

%
%

\subsection{Berry-Esseen Inequality}

If $F$ admits a Gibbs-Markov induced map of base $\Delta_0$ and return time function $R$, then for any $\varphi:M\to\mathbb R$ define $\varphi_{\Delta_0}:\Delta_0\to\mathbb R$ by
\[
\varphi_{\Delta_0}(x)=\sum_{i=0}^{ R(x)-1}\varphi (F^ix).
\]
%
%

Let $\varphi\in\mathcal H_\gamma$ be such that $\int \varphi d\mu=0$ and $\sigma^2$ be as in \eqref{eq.sigma}. Assume that $\sigma^2>0$ and that there exists $0<\delta\leq 1$ such that $\int |\varphi_{\Delta_0}|^2\chi_{|\varphi_{\Delta_0}|> z}d\mu\lesssim z^{-\delta}$, for large $z$. If $\delta=1$, assume also that  $\int |\varphi_{\Delta_0}|^3\chi_{|\varphi_{\Delta_0}|\leq z}d\mu$ is bounded. We say that Berry-Esseen Inequality holds for $\varphi$ if there exists $C>0$ such that for all $n\in\mathbb N$ and $a\in\mathbb R$ we have
\[
\left|\mu\left(\left\{x: \frac1{\sqrt n} \sum_{i=0}^{n-1} \varphi\circ F^i(x)\leq a\right\}\right) - \int_{-\infty}^a \frac1{\sigma\sqrt{2\pi}}\text e^{-\frac{x^2}{2\sigma^2}}dx\right|\leq \frac{C}{n^{\delta/2}}.
\]

%

\subsection{Almost Sure Invariance Principle}

Given $d\ge 1$ and a H\"older continuous $\varphi\colon M\to \mathbb R^d$ with mean zero, we
denote $$S_n=\sum_{i=0}^{n-1}\varphi\circ F^i, \quad\text{for each $n\ge 1$.}$$
We say that $\varphi$ satisfies an \emph{Almost Sure Invariance Principle} 
if there exists $\lambda > 0$ and  a probability
space supporting a sequence of random variables $\{S^*_n\}_n$ (which can be $\{S_n\}_n$ in the $d=1$ case)
and a $d$-dimensional Brownian motion
$W(t)$ such that
\begin{enumerate}
\item  $\{S_n\}_n $ and $ \{S^*_n\}_n $ are equally distributed;
\item  $S^*_n = W(n) + O(n^{1/2-\lambda})$, as $n\to\infty$,
almost everywhere.
\end{enumerate}

%
%
%

Satisfying an ASIP is a strong statistical property that implies other limiting laws such as the Central Limit Theorem, the Functional Central Limit Theorem or the Law of the Iterated Logarithm.
%

\end{document}